\documentclass[12pt,a4paper,reqno]{amsart}
\usepackage{amssymb,amscd}
\makeatletter
\@addtoreset{equation}{section}
\renewcommand{\@cite}[2]{{\rm[{{\textbf{#1}}\if@tempswa , #2\fi}]}}

\renewcommand{\@evenhead}{\thepage \hfil {\it Grigory Garkusha}}
\renewcommand{\@oddhead}{{\it FP-injective and weakly quasi-Frobenius rings} \hfil \thepage}
\renewcommand{\phi}{\varphi}
\renewcommand{\epsilon}{\varepsilon}

\renewcommand{\kappa}{\varkappa}
\DeclareMathOperator{\Hom}{Hom}
\DeclareMathOperator{\rad}{rad}
\DeclareMathOperator{\soc}{soc}
\DeclareMathOperator{\Ab}{Ab}
\DeclareMathOperator{\fg}{fg}
\DeclareMathOperator{\fp}{fp}
\DeclareMathOperator{\coh}{coh}
\DeclareMathOperator{\kr}{Ker}
\DeclareMathOperator{\Ext}{Ext}
\DeclareMathOperator{\im}{Im}

\DeclareMathOperator{\Mod}{Mod}
\DeclareMathOperator{\modd}{mod}
\newcommand{\f}{\mathfrak}
\newcommand{\lp}{\varinjlim}
\newcommand{\qq}[2]{#1_{\cc #2}}
\newcommand{\rg}{R(G)}
\newcommand{\rh}{R(H)}
\newcommand{\lra}[1]{\bl{#1}\longrightarrow\relax}
\newcommand{\bl}[1]{\buildrel #1\over}
\newcommand{\ii}[3]{0\longrightarrow#1\longrightarrow#2\longrightarrow#3}
\newcommand{\pp}[3]{#1\longrightarrow#2\longrightarrow#3\longrightarrow 0}
\newcommand{\cc}{\mathcal}

\newcommand{\lc}{{}_R{\cc C}}
\newcommand{\fl}{{\cc S}_R}
\newcommand{\ap}{{\cc S}^R}
\newcommand{\lfl}[1]{#1\in{\cc S}_R}
\newcommand{\lap}[1]{#1\in{\cc S}^R}
\newcommand{\iso}{\thickapprox}
\newcommand{\lcoh}[1]{#1\in\coh\lc}
\newcommand{\ten}[2]{{#1\otimes_R #2}}
\newcommand{\tenh}[2]{{#1\otimes_{\rh} #2}}
\newcommand{\tn}[2]{{#1\otimes #2}}
\newcommand{\tnn}[3]{{#1\bigotimes\limits_{#2} #3}}
\newcommand{\wh}{\widehat}
\newcommand{\wt}{\widetilde}
\newcommand{\les}[3]{0\longrightarrow#1\longrightarrow#2\longrightarrow#3\longrightarrow 0}

\newcommand{\rfp}{\modd-R}
\newcommand{\lfp}{R-\modd}
\newcommand{\Lfp}{R-\Mod}
\newcommand{\Rfp}{\Mod-R}
\newcommand{\homa}[2]{\Hom_R(#1,#2)}
\newcommand{\homg}[2]{\Hom_{\rg}(#1,#2)}
\newcommand{\homh}[2]{\Hom_{\rh}(#1,#2)}
\newtheorem{thm}{Theorem}[section]
\newtheorem{prop}[thm]{Proposition}
\newtheorem{cor}[thm]{Corollary}
\newtheorem{lem}[thm]{Lemma}

\newtheorem*{fpf}{IF-problem}
\newtheorem*{fgf}{FGF-problem}
\newtheorem*{cf}{CF-problem}
\newtheorem*{exs}{Examples}

\makeatother

\begin{document}

\footskip30pt

\setcounter{section}{-1}
\title{FP-injective and weakly quasi-Frobenius rings}
\author{Grigory Garkusha}
\address{Stary Peterhof, Russia}
\email{ggarkusha@mail.ru}
\begin{abstract}
The classes of $FP$-injective and weakly quasi-Frobenius rings are
investi\-gated. The properties for both classes of rings are closely linked
with embed\-ding of finitely presented mo\-dules in
$fp$-flat and free modules respectively. Using these
properties, we characterize the classes of coherent CF and
FGF-rings. Moreover, it is proved that the group ring
$\rg$ is $FP$-injective (weakly quasi-Frobenius) if
and only if the ring $R$ is $FP$-injective (weakly
quasi-Frobenius) and the group $G$ is locally finite.
\end{abstract}
\maketitle

\thispagestyle{empty}
\section{\it Introduction}

An application of the duality context with respect to the bimodule
${}_RR_R$ to the categories
of finitely generated left and right $R$-modules leads
to the case when $R$ is a noetherian self-injective ring. Such
rings are called quasi-Frobenius (or QF-rings). In turn,
an $R$-duality for categories of finitely presented modules
leads to the class of weakly quasi-Frobenius rings (or WQF-rings).
Such rings can be described as coherent $FP$-injective rings~\cite{GG}.

In the present paper we continue an investigation of the classes of
$FP$-injective and WQF-rings. To begin with, one must introduce
a notion of an $FP$-cogene\-rator, which plays an essential role in our
analysis, aproximately the same one as the notion of the cogenerator
for the class of QF-rings. Using also properties of $fp$-flat and
$fp$-injective modules, we give new criteria for both classes of rings
(theorems~\ref{www},~\ref{vvv}, and~\ref{uh}), which allow
to describe also the classes of coherent CF and FGF-rings.
Moreover, it is proved analogs of Renault's and Connell's
theorems for the $FP$-injective and weakly quasi-Frobenius group rings
respectively (theorems~\ref{gga} and~\ref{yes}).

It should be emphasized that the most difficult with the technical
point of view statements for the $FP$-injective rings are proved with
the help of the category of generalized $R$-modules
$\lc=(\rfp,\Ab)$ which consist of additive covariant functors
from the category of finitely presented right $R$-modules $\rfp$
to the category of abelian groups Ab. In our situation this is the
typical case since it is localizing subcategories
of the category $\lc$ and corresponding to them torsion functors
enable to adapt many properties we are interested in of the category of
modules to the
category of finitely presented modules. It is with the latter category
the most interesting statements for $FP$-injective and
WQF-rings are linked.

Throughout the paper the category of left (respectively right)
$R$-modules is denoted by $\Lfp$ (respectively $\Rfp$) and the category
of finitely presented left (respectively right) $R$-modules by $\lfp$
(respectively $\rfp$). The dual module $\homa MR$ of $M\in\Lfp$ is
denoted by $M^*$. Regular rings are supposed to be von Neumann regular.

I should like to thank A.~I.~Generalov for some helpful discussions.

\section{\it Preliminaries}

Recall that the {\em category of generalized left $R$-modules}
   $$\lc=(\rfp,\Ab)$$
consist of additive covariant functors from the category
of the finitely presented right $R$-modules $\rfp$ to the
category of abelian groups Ab. In this section we give
some properties of the category $\lc$ used later.
For more detailed information about the category
$\lc$ we refer the reader to~\cite{H} and here we, for
the most part, shall adhere to this paper.
All subcategories considered are supposed to be full.

We say that a subcategory $\cc S$ of an abelian category
$\cc C$ is a {\em Serre subcategory\/} if for every short
exact sequence
   $$\les XYZ$$
in $\cc C$ the object $Y\in\cc S$ if and only if $X$, $Z\in\cc S$.
A Serre subcategory $\cc S$ of a Grothendieck category
$\cc C$ is {\em localizing\/} if it is closed under taking
direct limits. Equivalently, the inclusion functor
$i:\cc S\to\cc C$ admits a right adjoint
$t=\qq tS:\cc C\to\cc S$ which takes every object
$X\in\cc C$ to the maximal subobject $t(X)$ of $X$ belonging to
$\cc S$. The functor $t$ one calls the {\em torsion functor}.

An object $X$ of a Grothendieck category $\cc C$ is {\em finitely generated\/}
if whenever there are subobjects $X_i\subseteq X$ with $i\in I$
satisfying $X=\sum_{i\in I}X_i$, then there is
a finite subset $J\subset I$ such that $X=\sum_{i\in J}X_i$.
The subcategory of finitely generated objects is denoted by
$\fg\cc C$. A finitely generated object $X$ is called
{\em finitely presented\/} if every epimorphism
$\gamma:Y\to X$ with $Y\in\fg\cc C$ has the finitely generaed kernel
$\kr\gamma$. By $\fp\cc C$ we denote the subcategory consisting of
finitely presented objects. Finally, we refer to a finitely presented object
$X\in\cc C$ as {\em coherent\/} if every finitely generated subobject
of $X$ is finitely presented. The corresponding subcategory of coherent
objects will be denoted by $\coh\cc C$.

The category $\lc$ is a locally coherent Grothendieck category, that is
every object $C\in\lc$ is a direct limit $C=\lp_IC_i$ of coherent
objects $\lcoh{C_i}$. Equivalently, the category $\coh\lc$ is abelian.
Moreover, $\lc$ has enough coherent projective generators
$\{(M,-)\}_{M\in\rfp}$. Thus, every coherent object $C\in\coh\lc$ has
a projective presentation
   $$\pp {(N,-)}{(M,-)}C,$$
where $M$, $N\in\rfp$.

We say that $M\in\lc$ is a {\em $\coh\lc$-injective object\/} if
$\Ext^1_{\lc}(C,M)=0$ for every $\lcoh C$. The fully faithful functor
$\ten-?:\Lfp\to\lc$, $M\mapsto\ten-M$, identifies the module category
$\Lfp$ with the subcategory of $\coh\lc$-injective objects of the
category $\lc$. In addition, the functor $\ten-M\in\coh\lc$
if and only if $M\in\lfp$. Furthermore, for every $\lcoh C$ there is
also an exact sequence
   $$\ii C{\ten-M}{\ten-N}$$
in $\coh\lc$ with $M$, $N\in\lfp$.

\begin{prop}[\cite{C,G}]\label{koh}
For a ring $R$ the following are equivalent:

$(1)$ $R$ is left coherent;

$(2)$ for every finitely presented right $R$-module $M$ the left
      $R$-module $M^*=\homa MR$ is finitely presented;

$(3)$ for every finitely presented right $R$-module $M$ the left
      $R$-module $M^*=\homa MR$ is finitely generated;

$(4)$ for every coherent object $\lcoh C$ the left $R$-module
      $C(R)$ is finitely presented;

$(5)$ for every coherent object $\lcoh C$ the left $R$-module
      $C(R)$ is finitely generated.
\end{prop}

Recall also that a monomorphism $\mu:M\to N$ in $\Lfp$ is a
{\em pure monomorphim\/} if for every $K\in\Rfp$ the morphism
$\tn K\mu$ is a monomor\-phism. Equivalently, the
$\lc$-morphism $\tn-\mu$ is a monomorphism.

In the sequel, we use the following Serre subcategories
of the category $\coh\lc$:
   \begin{gather*}
    {\cc S}^R=\{\lcoh C\mid C(R)=0\}\\
    {\cc S}_R=\{\lcoh C\mid (C,\ten-R)=0\},
   \end{gather*}
as well as the localizing subcategories $\vec{\cc S}^R$ and
$\vec{\cc S}_R$ of $\lc$
   \begin{gather*}
    \vec{\cc S}^R=\{C\in\lc\mid C=\lp C_i, C_i\in{\cc S}^R\}\\
    \vec{\cc S}_R=\{C\in\lc\mid C=\lp C_i, C_i\in{\cc S}_R\}.
   \end{gather*}
The corresponding $\vec{\cc S}^R$-torsion and
$\vec{\cc S}_R$-torsion functors will be denoted by
$t_{{\cc S}^R}$ and $t_{{\cc S}_R}$.

\section{\it $FP$-injective and weakly quasi-Frobenius rings}

A left $R$-module $M$ is said to be {\em $FP$-injective\/}
(or {\em absolutely pure\/}) if for every  $F\in\lfp$
we have: $\Ext^1_R(F,M)=0$, or equivalently,
every monomorphism $\mu:M\to N$ is pure~\cite[2.6]{St2}.
The ring $R$ is {\em left $FP$-injective\/} if the module
${}_RR$ is $FP$-injective. $M$ is an {\em $fp$-injective module\/}
if for every monomorphism $\mu:K\to L$ in $\lfp$ the morphism
$(\mu,M)$ is an epimorphism. Clearly, $FP$-injective modules are
$fp$-injective and every finitely presented $fp$-injective module
is $FP$-injective. $M$ is called {\em $fp$-flat\/} if for every
monomorphism $\mu:K\to L$ in $\rfp$ the morphism $\tn{\mu}M$ is
a monomorphism.

We refer to a left $R$-module $K$ as an
{\em $FP$-cogene\-rator\/} if for every non-zero homomorphism $f:M\to N$
from the finitely generated module $M$ to the finitely presented
module $N$ there exists $g\in\homa NK$ such that $gf\ne 0$. $K$ is said to
be an {\em $fp$-cogenerator\/} if for every non-zero homomorphism
$f:M\to N$ in $\lfp$
there exists $g\in\homa NK$ such that $gf\ne 0$. Obviosly,
$FP$-cogenerators are $fp$-cogenerators. On the other hand, it is not
hard to see that any $fp$-cogenerator is an $FP$-cogenerator when
the ring $R$ is left coherent.

\begin{lem}\label{xxx}
For a left $R$-module $K$ the following assertions are equivalent:

$(1)$ $K$ is an $FP$-cogenerator;

$(2)$ every finitely presented left $R$-module embeds in a product
      $K^I=\prod_IK$ of copies of the module $K$;

$(3)$ for every finitely presented left $R$-module $M$
      the following relation holds:
   $$\bigcap_{\phi\in\homa MK}\kr\phi=0.$$
\end{lem}

\begin{proof}
$(1)\Rightarrow (2)$. Let $M\in\lfp$, $I=\homa MK$, and let $\mu:M\to K^I$
be the homomorphism
such that $\mu(x)=(\phi(x))_{\phi\in I}$. If $0\ne x\in M$ and
$i:Rx\to M$ is an inclusion, there is $\phi:M\to K$
such that $\phi i(x)\ne 0$. Thus, $\mu$ is a monomorphism.

$(2)\Rightarrow (1)$. Let $0\ne f:M\to N$ be a homomorphism
from the finitely generated module $M$ to the finitely presented
module $N$. By assumption, there exists a monomophism
$g=(g_i)_{i\in I}:N\to K^I$. Then $gf\ne 0$, and so there is
$i_0\in I$ such that $g_{i_0}f\ne 0$.

$(1)\Rightarrow (3)$. Suppose $\mu:M\to K^I$ is the monomorphism
constructed in the proof of the implication $(1)\Rightarrow (2)$. Then
   $$\bigcap_{\phi\in\homa MK}\kr\phi=\kr\mu=0.$$
$(3)\Rightarrow (2)$. We may take $I=\homa MK$.
\end{proof}

Recall that a module $M\in\Lfp$ is {\em semireflexive\/}
(respectively {\em reflexive\/}) if the canonical homomorphism
$M\to M^{**}$ is a monomorphism (respectively isomorphism).

We are now in possession of all the information for proving
the following statement (cf.~\cite[2.5]{GG}):

\begin{thm}\label{www}
For a ring $R$ the following conditions are equivalent:

$(1)$ the module $R_R$ is $FP$-injective;

$(2)$ the module ${}_RR$ is an $FP$-cogenerator;

$(3)$ in $\Lfp$ there is an $fp$-flat $FP$-cogenerator;

$(4)$ in $\Lfp$ there is an $fp$-flat cogenerator;

$(5)$ if $\alpha:M\to N$ is a morphism in $\lfp$ such that
      $\alpha^*=\homa\alpha R$ is an epimorphism, then $\alpha$ is a monomorphism;

$(6)$ every finitely presented left $R$-module is semireflexive;

$(7)$ every finitely presented left $R$-module embeds in
      an $fp$-flat module;

$(8)$ every (injective) left $R$-module embeds in
      an $fp$-flat module;

$(9)$ every $FP$-injective left $R$-module is $fp$-flat;

$(10)$ every injective left $R$-module is $fp$-flat;

$(11)$ every indecomposable injective left $R$-module is $fp$-flat;

$(12)$ every flat right $R$-module is $fp$-injective.
\end{thm}

\begin{proof}
$(1)\Rightarrow (2)$. According to~\cite[2.5]{GG} $\fl\subseteq\ap$
in $\lc$. Consider a non-zero homomorphism $f:M\to N$ with $M\in\fg(\Lfp)$
and $N\in\lfp$. Suppose $C=\im(\tn-f)$; then $C$ is a finitely generated
subobject of the coherent object $\ten-N$. Therefore $\lcoh C$.
Assume that $gf=0$ for every $g\in N^*$. Consider an arbitrary
$\lc$-morphism $\gamma:C\to\ten-R$. Since $\ten-R$ is a
$\coh\lc$-injective object, there exists $\tn-h:\ten-N\to\ten-R$
such that $\tn-h|_C=\gamma$. But $hf=0$, and so $\gamma=0$.
Whence we obtain that $\lfl C$, and thus $\lap C$. We see that
$C(R)=0$, which yields $f=0$, a contradiction.

$(2)\Rightarrow (3)$. By assumption, the module ${}_RR$ is an
$FP$-cogenerator.

$(3)\Rightarrow (7)$. Since a direct product of $fp$-flat modules
is an $fp$-flat module (see~\cite[2.3]{GG}), our statement
follows from lemma~\ref{xxx}.

$(7)\Rightarrow (1)$. Let $f:M\to F$ be an embedding of a module $M\in\lfp$
in an $fp$-flat module $F$. Denote by $X=\kr(\tn-f)$. Let
$X=\sum_IX_i$, where each $X_i\in\fg\lc$. Since $X_i$ is a subobject
of the coherent object $\ten-M$, the object $X_i$ is coherent itself.
Because $X(R)=0$, each $X_i\in\ap$. Consequently, $X\in\vec\ap$.

Let us apply now the left exact $\fl$-torsion functor $t_{\fl}$
to the exact sequence
   $$\ii X{\ten-M}{\ten-F}.$$
Since $t_{\fl}(\ten-F)=0$ (see~\cite[2.2]{GG}), one gets that
$t_{\fl}(\ten-M)=t_{\fl}(X)\subseteq X\in\vec\ap$. Now let
$\lfl C$; since $C$ is a subobject of $\ten-M$ for some $M\in\lfp$,
from the relation
   $$C=t_{\fl}(C)\subseteq t_{\fl}(\ten-M)\in\vec\ap$$
we deduce that $\fl\subseteq\ap$, whence the ring $R$ is right
$FP$-injective by~\cite[2.5]{GG}.

$(1)\Rightarrow (5)$. Let $C=\kr(\tn-\alpha)$; then $\lcoh C$ and
since $\ten-R$ is a $\coh\lc$-injective object, there is an exact
sequence of abelian groups
   \begin{equation}\label{zq}
     (\ten-N,\ten-R)\lra{(\tn-\alpha,\ten-R)}
     (\ten-M,\ten-R)\lra{}(C,\ten-R)\lra{} 0.
   \end{equation}
Because $\alpha^*$ is an epimorphism, we conclude that
$\lfl C\subseteq\ap$. Thus $C(R)=0$, and hence
$\alpha$ is a monomorphism.

$(5)\Rightarrow (1)$. Let $\lfl C$; then there is an exact sequence
   $$0\to C\to\ten-M\bl{\tn-\alpha}\to\ten-N,$$
which induces an exact sequence of the form~\eqref{zq}.
Since $(C,\ten-R)=0$, $\alpha^*$ is an epimorphism, and hence
a monomorphism. So $\lap C$. By~\cite[2.5]{GG} the module $R_R$
is $FP$-injective.

$(1)\Leftrightarrow (6)$. This follows from~\cite[2.3]{Ja}.

$(4)\Rightarrow (3)$. Obvious.

$(1)\Rightarrow (9),(10),(11)$. This is a consequence of~\cite[2.5]{GG}.

$(9)\Rightarrow (10)\Rightarrow (11)$. Straightforward.

$(11)\Rightarrow (4)$. Suppose that $\{M_i\mid i\in I\}$ is a
system of representatives for isomorphism classes of simple
$R$-modules. Then every $E_i=E(M_i)$ is an indecomposable
injective $R$-module and the module $E=\prod_IE(M_i)$ is a
cogenerator in $\Lfp$. Because every module $E_i$ is $fp$-flat,
the module $E$ is $fp$-flat by~\cite[2.3]{GG}.

$(8)\Rightarrow (7)$. Easy.

$(10)\Rightarrow (8)$. It suffices to observe that the injective hull
$E(M)$ of the module $M$ is an $fp$-flat module.

$(1)\Rightarrow (12)$. This is a consequence of~\cite[2.5]{GG}.

$(12)\Rightarrow (1)$. Since the module $R_R$ is flat, it is
$fp$-injective, and so is $FP$-injective.
\end{proof}

A left (respectively right) ideal $I$ of the ring $R$ is {\em annulet\/}
if $I=\f l(X)$ (respectively $I=\f r(X)$), where $X$ is some subset
of the ring  $R$ and $\f l(X)=\{r\in R\mid rX=0\}$ (respectively
$\f r(X)=\{r\in R\mid Xr=0\}$). According to~\cite{F}
the left ideal $I$ is annulet if and only if $I=\f l\f r(I)$.

\begin{prop}\label{ok}
For a ring $R$ the following assertions hold:

\begin{itemize}
\item[$(1)$] if $R_R$ is an $FP$-injective module, then
\begin{itemize}
\item[(a)] for arbitrary finitely generated right ideals
$I$, $J$ of the ring $R$ one has: $\f l(I\cap J)=\f l(I)+\f l(J)$;

\item[(b)] for an arbitrary finitely generated left ideal
$I$ one has: $I=\f l\f r(I)$.
\end{itemize}

\item[$(2)$] if for $R$ the conditions\/ {\rm (a)} and\/ {\rm (b)}
from\/ $(1)$ hold, then every homomorphism of a finitely
generated right ideal of the ring $R$ into $R$ can be extended
to $R$.
\end{itemize}
\end{prop}

\begin{proof}
The proof is similar to that of~\cite[12.4.2]{K}.
\end{proof}

\begin{cor}\label{kkk}
Under the conditions (a) and (b) of proposition~\ref{ok}
if the ring $R$ is right coherent, then it is a right
$FP$-injective ring.
\end{cor}

\begin{proof}
By proposition~\ref{ok} we have $\Ext^1_R(R/I,R)=0$ for
every finitely generated right ideal $I$ of the ring $R$.
From~\cite[3.1]{St2} we deduce that $\Ext^1_R(M,R)=0$ for every
$M\in\rfp$, i.e.\ $R$ is right $FP$-injective.
\end{proof}

Now let us consider the situation when every finitely presented left
module embeds in a free module. Rings with such a property
one calls {\em left $IF$-rings.} In view of theorem~\ref{www} any left
$IF$-ring is right $FP$-injective. The next statement extends the list
of properties characterizing the $IF$-rings (cf.~\cite{C}).

\begin{prop}\label{qqq}
For a ring $R$ the following conditions are equivalent:

$(1)$ $R$ is a left $IF$-ring;

$(2)$ every finitely presented left $R$-module embeds in a flat
      $R$-module;

$(3)$ every $FP$-injective left $R$-module is flat;

$(4)$ every injective left $R$-module is flat.
\end{prop}

\begin{proof}
$(1)\Rightarrow (2)$ is trivial.

$(2)\Rightarrow (1)$. Let $f:K\to M$ be an embedding of a module
$K\in\lfp$ in a flat module $M$. By theorem of Govorov and
Lazard~\cite[11.32]{F}
the module $M$ is a direct limit $\lp_IP_i$ of the projective modules
$P_i$. By~\cite[V.3.4]{St} there is $i_0\in I$ such that $f$
factors through $P_{i_0}$. Therefore $K$ is a submodule of
$P_{i_0}$. It remains to observe that $P_{i_0}$ is a submodule
of some free module.

$(1)\Leftrightarrow (4)$. This follows from~\cite[2.1]{C}.

$(3)\Rightarrow (4)$. Obvious.

$(4)\Rightarrow (3)$. If $M$ is an $FP$-injective left $R$-module,
then the sequence
   $$\les ME{E/M},$$
in which $E=E(M)$, is pure. Since $E$ is a flat module, the
module $E/M$ is flat by~\cite[I.11.1]{St}. Let
$\wh M=\Hom_{\Bbb Z}(M,\Bbb Q/\Bbb Z)$ denote the character module of
$M$. By~\cite[I.10.5]{St} the modules $\wh E$ and $\wh{E/M}$
are injective, and so the exact sequence
   $$\les{\wh {E/M}}{\wh E}{\wh M}$$
splits. Consequently, the module $\wh M$ is injective,
and so $M$ is flat by~\cite[I.10.5]{St}.
\end{proof}

Colby~\cite{C} has constructed an example of a left $IF$-ring,
which is not a right $IF$-ring.

\begin{prop}\label{do}
If $R$ is a left $FP$-injective ring and a left $IF$-ring, then
it is right coherent.
\end{prop}

\begin{proof}
By assumption, every $K\in\lfp$ embeds in a free module
(and so in a finitely generated free module as well).
One has the following exact sequence
   \begin{equation}\label{2.1}
    \ii K{R^n}{R^m}
   \end{equation}
in $\lfp$. Since the module ${}_RR$ is $FP$-injective, one gets
an exact sequence
   \begin{equation}\label{2.2}
    \pp {R^m}{R^n}{K^*}
   \end{equation}
in $\rfp$, hence $K^*\in\rfp$. By proposition~\ref{koh} the ring $R$
is right coherent.
\end{proof}

\begin{fpf}{\rm
Is it true that any left $IF$-ring is right coherent?
}\end{fpf}

It should be remarked that $IF$-problem, in view of proposition~\ref{qqq},
is equivalent to Jain's problem~\cite[p.~442]{Ja}:
will be the ring $R$ right coherent if every injective
left $R$-module is flat?

We recall that the ring $R$ is {\em almost regular\/} (see~\cite{GG})
if every (both left and right) module is $fp$-flat.
By theorem~\ref{www} almost regular rings are two-sided
$FP$-injective rings.

\begin{cor}
An almost regular ring $R$ will be a left $IF$-ring if and only
if it is regular.
\end{cor}

\begin{proof}
Clearly, an almost regular ring is regular if and only if
it is left or right coherent. Therefore our assertion immediately
follows from proposition~\ref{do}.
\end{proof}

Recall also that the ring $R$ is {\em indiscrete\/} if it is
a simple almost
regular ring. Prest, Rothmaler and Ziegler~\cite{P} have constructed
an example of a non-regular indiscrete ring.

The ring $R$ is said to be {\em weakly quasi-Frobenius\/} (or
{\em {\rm WQF}-ring\/}) if it determines an $R$-duality between the
categories of finitely presented left and right $R$-modules. Such
rings can be described as (left and right) $FP$-injective
(left and right) coherent rings~\cite[2.11]{GG}.

The next theorem extends the list of properties characterizing
the WQF-rings (cf.~\cite[2.11; 2.12]{GG}, \cite[2.2]{C}).

\begin{thm}\label{vvv}
For a ring $R$ the following assertions are equivalent:

$(1)$ $R$ is a {\rm WQF}-ring;

$(2)$ $R$ is a left and right $IF$-ring;

$(3)$ $R$ is left coherent and left $FP$-injective, and every
      cyclic finitely presented left $R$-module embeds in
      a free module;

$(4)$ every left and every right $FP$-injective $R$-module is flat;

$(5)$ every left and every right injective $R$-module is flat.
\end{thm}

\begin{proof}
$(1)\Rightarrow (4)$. This follows from~\cite[2.11]{GG}.

$(2)\Leftrightarrow (4)\Leftrightarrow (5)$. Apply
proposition~\ref{qqq}.

$(2)\Rightarrow (1)$. Since any left and right $IF$-ring is a two-sided
$FP$-injective ring, our statement follows from proposition~\ref{do}.

$(1)\Rightarrow (3)$. Straightforward.

$(3)\Rightarrow (1)$. Since for every cyclic $K\in\lfp$
the dual module $K^*\ne 0$, the proof of right $FP$-injectivity
of the ring $R$ is similar to that of~\cite[2.9]{GG}.

Let us show that the ring $R$ is right coherent. In view of
proposition~\ref{koh} it suffices to prove that $K^*\in\rfp$
for every $K\in\lfp$. We use induction on the number of generators
$n$ of the module $K$. When $n=1$, considerating exact sequences~\eqref{2.1}
and~\eqref{2.2} for $K$, one gets $K^*\in\rfp$. If $K$ is finitely
presented on $n$ generators, let $K'$ be the submodule of $K$ generated
by one of these generators. Since $R$ is left coherent, the modules
$K'$ and $K/K'$ are finitely presented on less than $n$ generators.
Because $R$ is left $FP$-injective, one has an exact sequence
   $$\les{(K/K')^*}{K^*}{(K')^*},$$
where both $(K')^*$ and $(K/K')^*$ are finitely presented by induction.
Thus $K^*\in\rfp$.
\end{proof}

Now we combine the preceding arguments in the following theorem (cf.\
properties of QF-rings~\cite[24.4]{F2}, \cite[13.2.1]{K} and also
properties of rings with full duality~\cite[12.1.1]{K}):

\begin{thm}\label{uh}
For a two-sided coherent ring $R$ the following conditions are equivalent:

$(1)$ $R$ is a {\rm WQF}-ring;

$(2)$ the modules ${}_RR$ and $R_R$ are $FP$-injective;

$(3)$ ${}_RR$ and $R_R$ are $FP$-cogenerators;

$(4)$ the module $R_R$ is an $FP$-injective $FP$-cogenerator;

$(5)$ every left and every right finitely presented $R$-module is reflexive;

$(6)$ every left and every right cyclic finitely presented $R$-module
      is reflexive;

$(7)$ every left and every right cyclic finitely presented $R$-module
      embeds in a free module;

$(8)$ for a finitely generated left ideal $I$ and for a finitely
      generated right ideal $J$ of the ring $R$ one has:
      $\f l\f r(I)=I$ and $\f r\f l(J)=J$.
\end{thm}

\begin{proof}
$(1)\Leftrightarrow (2)$. This follows from~\cite[2.11]{GG}.

$(2)\Leftrightarrow (3)\Leftrightarrow (4)$. Apply theorem~\ref{www}.

$(2)\Leftrightarrow (5)$. This follows from~\cite[4.9]{St2}.

$(5)\Rightarrow (6)$. Obvious.

$(6)\Rightarrow (7)$. Let $M$ be a cyclic finitely presented left
$R$-module. In view of proposition~\ref{koh} the module $M^*\in\rfp$,
and so there is an epimorphism $R^n\to M^*$. Hence
$M\iso M^{**}\to R^n$ is a monomorphism.

$(7)\Rightarrow (8)$. Let $I$ be a finitely generated left ideal of the
ring $R$. By assumption, the module $R/I$ embeds in a free module.
By~\cite[20.26]{F2} there exists a finite subset $X$ of $R$ such that
$I=\f l(X)$, i.~e., $I$ is an annulet ideal. By symmetry,
every finitely generated right ideal is annulet.

$(8)\Rightarrow (2)$. In view of corollary~\ref{kkk} it suffices
to show that for arbitrary finitely generated right ideals $I$ and
$J$ of the ring $R$ the following equality holds:
$\f l(I\cap J)=\f l(I)+\f l(J)$. Since $R$ is coherent by assumption,
by the Chase theorem~\cite[I.13.3]{St} both $I\cap J$ and
$\f l(I)+\f l(J)$ are finitely generated ideals.

One has
   $$\f r\f l(I\cap J)=I\cap J=\f r\f l(I)\cap\f r\f l(J)=
     \f r(\f l(I)+\f l(J)).$$
Applying $\f l$, one gets
   $$\f l(I\cap J)=\f l\f r(\f l(I)+\f l(J))=\f l(I)+\f l(J).$$
Thus $R_R$ is $FP$-injective. Likewise, ${}_RR$ is $FP$-injective.
\end{proof}

It is well-known that QF-rings have the global dimension to be equal to
0 (and then the ring $R$ is semisimple), or $\infty$. In turn,
WQF-rings, in view of~\cite[3.6]{St2}, have the weak global
dimension to be equal to 0 (and then the ring $R$ is regular), or $\infty$.

Some examples of WQF-rings the reader can find in~\cite{C}.

Next, we consider the class of rings over which every
finitely generated left $R$-module embeds in a free $R$-module.
Such rings we shall call {\em left\/ {\rm FGF}-rings.}
Clearly, any FGF-ring will be an $IF$-ring.
In turn, if every cyclic left $R$-module embeds in a free
$R$-module, the ring $R$ one calls a {\em left\/ {\rm 'F}-ring.}
The following problems are still open (see~\cite{RS}):

\begin{fgf}{\rm
Does the class of left FGF-rings coincide with the class of QF-rings?
}\end{fgf}

\begin{cf}{\rm
Will be left CF-rings left artinian?
}\end{cf}

Recall also that the ring $R$ is a {\em left Kasch ring\/} if the injective
hull $E({}_RR)$ of the module ${}_RR$ is an injective cogenerator
in $\Lfp$. Equivalently, for every non-zero cyclic left $R$-module
$M$ the dual module $M^*\ne 0$ (see~\cite[XI.5.1]{St}).

\begin{lem}\label{ttt}
For a ring $R$ the following assertions hold:

$(1)$ if $R$ is a left coherent ring and a left\/ {\rm CF}-ring,
      then it is a left noetherian ring and a left Kasch ring;

$(2)$ if $R$ is a left noetherian ring and a left $IF$-ring,
      then it is a left\/ {\rm FGF}-ring. If $R$
      is a left coherent and a left\/ {\rm FGF}-ring, then it is
      a left noetherian ring and a left $IF$-ring.
\end{lem}

\begin{proof}
(1). It is easy to see that $R$ is a left Kasch ring. Thus
we must show that the module ${}_RR$ is noetherian.

Suppose $I$ is a left ideal of the ring $R$. By assumption,
the module $R/I$ is a submodule of a free $R$-module $R^n$
for some $n\in\Bbb N$. Since the ring $R$ is left coherent,
the module $R^n$ is coherent, and hence the module $R/I$
is finitely presented, i.e., $I$ is a finitely generated ideal.

(2). It is necessary to observe that over a noetherian ring every
finitely generated module is finitely presented and also make use of
the first statement.
\end{proof}

\begin{prop}
For a two-sided coherent ring $R$ the following assertions
are equivalent:

$(1)$ $R$ is a left\/ {\rm FGF}-ring;

$(2)$ the module ${}_RR$ is a noetherian $FP$-cogenerator;

$(3)$ the module ${}_RR$ is noetherian and the module $R_R$
      is $FP$-injective;

$(4)$ $R$ is a left noetherian ring, a left Kasch ring, and
      the module $E({}_RR)$ is flat.
\end{prop}

\begin{proof}
$(1)\Rightarrow (2),(3)$. This follows from lemma~\ref{ttt}, lemma~\ref{xxx}
and theorem~\ref{www}.

$(1)\Rightarrow (4)$. Apply lemma~\ref{ttt} and proposition~\ref{qqq}.

$(2)\Leftrightarrow (3)$. This is a consequence of theorem~\ref{www}.

$(2)\Rightarrow (1)$. Since $R$ is a left noetherian ring,
every finitely generated left $R$-module is finitely presented.
By lemma~\ref{xxx} every finitely presented left $R$-module
is a submodule of the module $R^I$. Because
the ring $R$ is right coherent, the module $R^I$ is flat
by~\cite[I.13.3]{St}. By proposition~\ref{qqq} $R$ is a left $IF$-ring and
by lemma~\ref{ttt} $R$ is also a left FGF-ring.

$(4)\Rightarrow (1)$. In this case the proof is similar to the proof of
the implication $(2)\Rightarrow (1)$ if we observe that every finitely
generated left $R$-module is a submodule of the flat module $E^I$.
\end{proof}

The ring $R$ is called {\em left semiartinian\/} if every
non-zero cyclic left $R$-module has a non-zero socle. $R$
is {\em semiregular\/} if $R/\rad R$ is a regular ring.

\begin{lem}\cite[2.1]{RS}\label{gr}
A left\/ {\rm CF}-ring $R$ is left semiartinian if and only
if $\soc({}_RR)$ is an essential submodule in ${}_RR$.
\end{lem}

The next two statements partially solve CF- and FGF-problems
respectively (cf.~\cite{RS}):

\begin{prop}\label{well}
Suppose $R$ is a left noetherian ring and a left\/
{\rm CF}-ring. Then the following conditions are equivalent
for $R$:

$(1)$ $R$ is left artinian;

$(2)$ $R$ is left or right semiartinian;

$(3)$ $R$ is a semiregular ring;

$(4)$ $R$ is a semiperfect ring;

$(5)$ $\soc({}_RR)$ is an essential submodule in ${}_RR$.
\end{prop}

\begin{proof}
$(1)\Rightarrow (2)$. Any left artinian ring is left perfect,
and so is right semiartinian by~\cite[VIII.5.1]{St}.

$(2)\Rightarrow (1)$. Our assertion follows from~\cite[VIII.5.2]{St}.

$(1)\Rightarrow (3)$. Easy.

$(3)\Rightarrow (4)$. Since $\soc({}_RR)$ is a finitely
generated left ideal of $R$, our assertion follows from~\cite[4.1]{RS}.

$(4)\Rightarrow (5)$. Obvious.

$(5)\Rightarrow (1)$. By the preceding lemma $R$ is left
semiartinian and since $R$ is left noetherian,
our assertion follows from~\cite[VIII.5.2]{St}.
\end{proof}

\begin{prop}
Let $R$ be a left noetherian ring and a left\/ {\rm FGF}-ring.
Then the following are equivalent for $R$:

$(1)$ $R$ is a {\rm QF}-ring;

$(2)$ $R$ is a {\rm WQF}-ring;

$(3)$ $R$ is a left $FP$-injective ring;

$(4)$ $R$ is a right Kasch ring;

$(5)$ $R$ is a semiregular ring;

$(6)$ $R$ is left or right semiartinian;

$(7)$ $\soc({}_RR)$ is an essential submodule in ${}_RR$.
\end{prop}

\begin{proof}
The implications $(1)\Rightarrow (2)\Rightarrow (3)$,
$(1)\Rightarrow (4),(5)$ are obvious.

$(3)\Rightarrow (1)$. Our assertion follows from~\cite[2.6]{GG}.

$(4)\Rightarrow (3)$. Over a right Kasch ring the module
$M^*\ne 0$ for every non-zero cyclic right $R$-module $M$.
Since $R$ is a left FGF-ring, by theorem~\ref{www} $R$ is a right
$FP$-injective ring. Therefore $R$ is a left $FP$-injective ring
by~\cite[2.9]{GG}.

$(5)\Leftrightarrow (6)\Leftrightarrow (7)$. Apply proposition~\ref{well}.

$(7)\Rightarrow (1)$. By proposition~\ref{well} the ring $R$ is left
artinian. Since any left artinian ring is right
perfect, our assertion follows from~\cite[2.5]{RS}.
\end{proof}

\section{\it Group rings}

Let $R$ be a ring and $G$ a group. Denote the group ring
of $G$ with coefficients in $R$ by $\rg$.

\begin{lem}\label{mmm}
Suppose $G$ is a finite group and $M\in\Mod-R(G)$. Let
$f=(f_i)_I:M\to R^I$ is some $R$-monomorphism with $I$ some set of
indices. Then the $R(G)$-homomorphism $\wt f=(\wt f_i)_I:M\to R(G)^I$
defined by the rule
   $$\wt f_i(m)=\sum\nolimits_{g\in G}f_i(mg)g^{-1}$$
is an $R(G)$-monomorphism.
\end{lem}

\begin{proof}
It is directly verified that $\wt f$ is indeed a monomorphism
of $R(G)$-modules.
\end{proof}

According to Renault's theorem~\cite{R} the group ring
$\rg$ is left self-injective if and only if the ring
$R$ is left self-injective and the group $G$ is finite.
In turn, for the $FP$-injective rings there is the following
statement.

\begin{thm}\label{gga}
The group ring $\rg$ is left $FP$-injective if and only if
the ring $R$ is left $FP$-injective and the group $G$ is locally finite.
\end{thm}

\begin{proof}
Suppose $\rg$ is left $FP$-injective and $M\in\rfp$.
Then the module $\ten M{\rg}\in\modd-\rg$ and, in view of lemma~\ref{xxx}
and theorem~\ref{www}, there exists a monomorphism
$\mu:\ten M{\rg}\to\rg^I$ with $I$ some set of indices.
Then the composition
   $$M\lra\beta\ten M{\rg}\lra\mu\rg^I$$
of morphisms $\mu$ and $\beta$ is an $R$-monomorphism, where
$\beta(m)=\tn me$, $e$ is a unit of the group $G$. Since $\rg$ is a
free $R$-module, the $R$-module $\rg^I$ is $fp$-flat by~\cite[2.3]{GG}.
Theorem~\ref{www} implies that the ring $R$ is left $FP$-injective.

Let us show now that the group $G$ is locally finite.
Let $H$ be a non-trivial subgroup of the group $G$ generated by
elements $\{h_i\}_{i=1}^n$. Then the right ideal $\omega H$ of the ring
$\rg$ generated by the elements $\{1-h_i\}_{i=1}^n$ is non-zero.
By proposition~\ref{ok} $\f l(\omega H)\ne 0$, and so $H$ is
finite by~\cite[2.1]{Con}.

Now let $R$ be left $FP$-injective and the group $G$ locally finite.
To begin, let us show that the ring $\rg$ is left $FP$-injective
if $G$ is finite.

Suppose $M\in\modd-\rg$. Because the group $G$ is finite, $M\in\rfp$.
Since the module $R_R$ is an $FP$-cogenerator by theorem~\ref{www},
by lemma~\ref{xxx} $M_R$ is a submodule of $R^I$, where $I$ is some
set. By lemma~\ref{mmm} $M_{\rg}$ is a submodule of $\rg^I$.
Consequently, $\rg$ is an $FP$-cogenerator, and hence
$\rg$ is a left $FP$-injective ring by theorem~\ref{www}.

Next, suppose that $G$ is an arbitrary locally finite group
and $M\in\rg-\modd$. Then there is a short exact sequence
of $\rg$-modules
   $$\les{K\bl i}{\rg^n\bl p}M.$$
Let $X$ be a finite set of generators for ${}_{\rg}K$. Because
$G$ is locally finite, there is a finite subgroup $H$ of
$G$ such that $\rh X\subseteq\rh^n\subseteq\rg^n$. We result
in the short exact sequence of $\rh$-modules
   \begin{equation}\label{en}
    \les{K'\bl{\bar i}}{\rh^n\bl{\bar p}}{M'},
   \end{equation}
where $K'=\rh X$, $M'=\rh Y$,
$Y$ is a finite set of generators of the module ${}_{\rg}M$.
Tensoring sequence~\eqref{en} on $\rg_{\rh}$, one gets
the following commutative diagram with exact rows:
   $$\begin{CD}
      0@>>>\tnn{\rg}{\rh}{K'}@>\tn 1{\bar i}>>\tnn{\rg}{\rh}{\rh^n}@>\tn 1{\bar p}>>\tnn{\rg}{\rh}{M'}@>>>0\\
      @.@V{\alpha}VV@VV{\beta}V@VV{\gamma}V\\
      0@>>>K@>>i>\rg^n@>>p>M@>>>0,\\
     \end{CD}$$
in which $\beta$ is an isomorphism, $\alpha(\tn rx)=rx$. Clearly, $\alpha$
is an isomorphism, and hence $\gamma$ is an isomorphism.

If we showed that every $f\in\homg K{\rg}$ is extended
to some $\sigma\in\homg{\rg^n}{\rg}$, we would obtain
that $\Ext^1_{\rg}(M,\rg)=0$, as required.

So suppose that $f\in\homg K{\rg}$ and $\tau\in\homh{\rg}{\rh}$,
$\tau(\sum_{g\in G}r(g)g)=\sum_{g\in H}r(g)g$. Consider
$\bar f=\tau f|_{K'}\in\homh{K'}{\rh}$. Because $\rh$ is a left $FP$-injective ring,
there is $\bar\sigma:\rh^n\to\rh$ such that $\bar f=\bar\sigma\bar i$.
One gets
   $$f=(\tn 1{\bar f})\alpha^{-1}=
     (\tn 1{\bar\sigma})(\tn 1{\bar i})\alpha^{-1}=
     (\tn 1{\bar\sigma})\beta^{-1}i=\sigma i,$$
where $\sigma=(\tn 1{\bar\sigma})\beta^{-1}$ is
the required homomorphism.
\end{proof}

The Renault theorem and theorem~\ref{gga} implies that given
an arbitrary self-injective ring $R$, one can construct
$FP$-injective rings which are not self-injective. To be definite,
the following statement holds:

\begin{cor}
If $R$ is a left self-injective ring, $G$ is a locally finite
group, and $|G|=\infty$, then the group ring $\rg$ is
left $FP$-injective but not left self-injective.
\end{cor}

\begin{prop}[Colby~\cite{C}]\label{bn}
The group ring $\rg$ is a left $IF$-ring if and only if
$R$ is a left $IF$-ring and the group $G$ is locally finite.
\end{prop}

\begin{proof}
If $\rg$ is a left $IF$-ring, it is also right $FP$-injective
by theorem~\ref{www}. By the preceding theorem the group
$G$ is locally finite. Similar to the proof of theorem~\ref{gga},
given $M\in\lfp$ there is a composition of $R$-monomorphisms
   $$M\lra\beta\ten{\rg}M\lra\mu\rg^n$$
with $\beta(m)=\tn em$. Since $\rg$ is a free $R$-module,
the module $M$ is a submodule of the free $R$-module $\rg^n$.

Conversely, let $R$ be a left $IF$-ring and $G$ a locally finite group.
First, let us prove that $\rg$ is a left $IF$-ring if
$G$ is a finite group. For this consider $M\in\rg-\modd$. Since
$G$ is finite, $M\in\lfp$, and so ${}_RM$ is a submodule of $R^n$ for
some $n\in\Bbb N$. By lemma~\ref{mmm} ${}_{\rg}M$ is a submodule of
$\rg^n$, i.¥., $\rg$ is indeed a left $IF$-ring.

If $G$ is an arbitrary locally finite group, then for any
$M\in\rg-\modd$ there is a finite subgroup $H$ of the group $G$ such that
   $$M\iso\tenh{\rg}{\rh Y}$$
with $Y$ a finite set of generators for ${}_{\rg}M$
(see the proof of theorem~\ref{gga}). By assumption,
the $\rh$-module $\rh Y$ is a submodule of a free module
$\rh^n$ for some $n\in\Bbb N$. Thus, $M$ is a submodule of the free module
$\rg^n\iso\tenh{\rg}{\rh^n}$.
\end{proof}

\begin{thm}\label{yes}
The group ring $\rg$ is weakly quasi-Frobenius if and only if
the ring $R$ is weakly quasi-Frobenius and the group $G$
is locally finite.
\end{thm}

\begin{proof}
Theorem~\ref{vvv} implies that any WQF-ring is a left and right
$IF$-ring. Therefore our statement immediately follows
from proposition~\ref{bn}.
\end{proof}

It is well-known that the group ring $\rg$ is semisimple
(see~\cite{Con,F}) if and only if the ring $R$ is semisimple,
the group $G$ is finite, and $|G|$ is invertible in $R$.
In turn, by theorem of Auslander and McLaughlin
(see~\cite{Con}) $\rg$ is a regular ring if and only if
the ring $R$ is regular, the group $G$ is locally finite,
and for every finite subgroup $H$ of $G$ the equality
$|H|=n$ implies $nR=R$.

To conclude, we give some examples of WQF-rings which are
simultaneously neither QF-rings, nor regular rings.

\begin{exs}{\rm
(1) Given an arbitrary regular ring $R$, we can construct WQF-rings
which will not be regular. Namely, it is necessary to consider
an arbitrary locally finite group $G$, in which there is at least one
finite subgroup $H$ of $G$ such that the order $|H|$ is not
a unit in $R$.

To take an example, consider the field $K$ of the characteristic
$p\ne 0$. Let $R=\prod_{i=1}^\infty K_i$, $K_i=K$, be the ring
with component-wise operations. Then $R$ is a regular
but not semisimple ring, as one easily sees.
If $G$ is a finite group such that $p$ devides $|G|$, then
the ring $\rg$ is a weakly quasi-Frobenius ring being
neither quasi-Frobenius, nor regular.

(2) Let $R$ be an arbitrary QF-ring, $G$ an arbitrary
locally finite group, and $|G|=\infty$. Then $\rg$ is
a weakly quasi-Frobenius ring but not quasi-Frobenius.
Moreover, $\rg$ is regular if and only if $R$ is a semisimple
ring and the order of every finite subgroup of $G$
is invertible in $R$.

As an example, if $K$ is the field of the characteristic
$p\ne 0$, the group $G=\cup_{k\ge 1}G_k$, where every
$G_k$ is a cyclic group with a generator $a_k$ of the order
$p^k$, and $a_k=a_{k+1}^p$, then the group algebra
$K(G)=\lp_kK(G_k)$ is weakly quasi-Frobenius (see also~\cite{M})
being neither quasi-Frobenius, nor regular.
}\end{exs}

\end{document}